\documentclass{amsart}

\usepackage{amsmath,amssymb,amscd,amsfonts}
\usepackage[all]{xy}

\newtheorem{theorem}{Theorem}
\newcommand{\bt}{\begin{theorem}}
\newcommand{\et}{\end{theorem}}
\newtheorem{lemma}{Lemma}
\newcommand{\bl}{\begin{lemma}}
\newcommand{\el}{\end{lemma}}
\newtheorem{corollary}{Corollary}
\newcommand{\bc}{\begin{corollary}}
\newcommand{\ec}{\end{corollary}}
\newcommand{\beq}{\begin{equation}}
\newcommand{\eeq}{\end{equation}}
\newcommand{\benum}{\begin{enumerate}}
\newcommand{\eenum}{\end{enumerate}}

\newcommand{\Z}{\ensuremath{\mathbf Z}}
\newcommand{\Q}{\ensuremath{\mathbf Q}}

\newcommand{\mbF}{\ensuremath{ \mathbf F}}

\newcommand{\mce}{\ensuremath{ \mathcal E}}

\DeclareMathOperator{\card}{\text{card}}

\newcommand{\bmat}{\left(\begin{matrix}}
\newcommand{\emat}{\end{matrix}\right)}

\title{An inverse problem for finite Sidon sets}

\author{Melvyn B. Nathanson}
\address{Lehman College (CUNY), Bronx, New York 10468}
\email{melvyn.nathanson@lehman.cuny.edu}



\thanks{Supported in part by a grant from the PSC-CUNY Research Award Program.}

\date{\today}

\begin{document}

\begin{abstract} 
Here is a  direct problem for Sidon sets: 
Given a linear form  $\varphi =  c_1 x_1 + \cdots + c_h x_h $, 
 construct and describe sets $A$ that are Sidon sets for $\varphi$.  
This paper considers an inverse problem for Sidon sets: 
Given a set $A$,  determine the linear forms $\varphi$ 
such that $A$ is a Sidon set for $\varphi$. 
\end{abstract}

\maketitle

\section{Sidon sets for linear forms}               \label{LinearPerturb:section:forms}
A classical Sidon set is a set $A$ of positive integers such that every integer has at most one representation 
as the sum of two elements of $A$.  Equivalently, $A$ has the property that if $a_1,a_2,a'_1,a'_2 \in A$ and 
$a_1+a_2 = a'_1 +a'_2$, then $\{a_1,a_2\} = \{a'_1,a'_2\}$. 
In this paper we consider an analogue of Sidon sets for  linear forms.  

Let \mbF\ be a field, let $\mbF^{\times} = \mbF\setminus \{ 0\}$,  
and let $h$ be a positive integer.  
We consider linear forms  
\[      
\varphi = \varphi(x_1, x_2,\ldots, x_h) =  c_1 x_1 + c_2 x_2 +\cdots + c_h x_h  
\]
where $c_i \in \mbF^{\times}$ for all $i \in \{1,\ldots, h\}$.  
The linear form $\varphi$ is \emph{monic} if $c_1 = 1$.
A linear form with 2 variables is \emph{binary}.  
A linear form with $h$ variables is called \emph{$h$-ary}.

Let  $A$ be a nonempty subset \mbF.  
For $c,t \in \mbF$, the \emph{$c$-dilate} of $A$ is the set 
\[
c\ast A = \{ca:a\in A\}
\]
and the \emph{$t$-translate} of $A$ by $t$ is the set 
\[
A+t= \{ a+t :a\in A\}.
\]
Let 
\[
A^h = \left\{ (a_1,\ldots, a_h): a_i \in A \text{ for all } i \in \{1,\ldots, h\}  \right\}
\]
be the set of all $h$-tuples of elements of $A$.  
The \emph{$\varphi$-image of $A$} is the set 
\begin{align*}
\varphi(A) & = \left\{ \varphi(a_{1},\ldots, a_{h}): (a_1,\ldots, a_h) \in A^h \right\} \\ 
& = \left\{ c_1 a_{1} + \cdots +  c_h a_{h}: (a_1,\ldots, a_h) \in A^h \right\} \\
& = c_1\ast A + \cdots + c_h \ast A. 
\end{align*}
Thus, $\varphi(A)$ is a sum of dilates.  
We define $\varphi(\emptyset) = \{0\}$.

A nonempty subset $A$  of $\mbF$ is a 
\emph{Sidon set for the linear form $\varphi$} or, simply, a 
\emph{$\varphi$-Sidon set} 
if it satisfies the following property:  
For all $h$-tuples $(a_1,\ldots, a_h) \in A^h$ and $ (a'_1,\ldots, a'_h) \in A^h$,   
if 
\[
\varphi(a_1,\ldots, a_h) = \varphi(a'_1,\ldots, a'_h) 
\]
then $(a_1,\ldots, a_h) = (a'_1,\ldots, a'_h)$, 
that is, $a_i = a'_i$ for all $i \in \{1,\ldots, h\} $.  
Thus, $A$ is a $\varphi$-Sidon set if the linear form $\varphi$ is one-to-one on $A^h$.  
If $A$ is a finite set with $k$ elements, then $A$ is a $\varphi$-Sidon set for the $h$-ary linear form 
$\varphi$ if and only if $|\varphi(A)| = k^h$.

Of special interest is the case $\mbF = \Q$ with $\varphi$-Sidon sets of integers or positive integers.

A \emph{direct problem} for Sidon sets for linear forms is the following: 
Given a linear form $\varphi$, describe the subsets of $\mbF$ that are $\varphi$-Sidon sets.
In this paper we consider an \emph{inverse problem}:    
Given a subset $A$ of $\mbF$, describe the linear forms $\varphi$ such that 
$A$ is a $\varphi$-Sidon set.  

Sidon sets for linear forms were introduced by Nathanson~\cite{nath07m,nath2021-200}.

\section{Monic forms and normalized sets} 

Associated with the linear form
\[
\varphi = c_1x_1+c_2x_2 + \cdots + c_hx_h 
\]
is the monic linear form 
\[
\psi = x_1+ c'_2x_2 + \cdots + c'_hx_h 
\]
where 
\[
c'_i = \frac{c_i}{c_1} \qquad \text{for $i = 2,3,\ldots, h$.}
\]
Let $A$ be a nonempty subset of $V$.  
For all $(a_1,\ldots, a_h) \in A^h$, we have 
\begin{align*}
\varphi(a_1,\ldots, a_h) 
& = c_1a_1+c_2a_2 + \cdots + c_ha_h \\
& = c_1\left( a_1+ \frac{c_2}{c_1}a_2 + \cdots +  \frac{c_h}{c_1} a_h\right) \\
& = c_1 \psi(a_1,\ldots, a_h) 
\end{align*}
and so 
\[
 \varphi(A) = c_1\ast \psi(A).
\]
Because $c_1 \neq 0$, the set $A$ is a  $\varphi$-Sidon set if and only if $A$ is a $\psi$-Sidon set. 
Thus, it suffices to consider only monic linear forms.

A subset $A$ of \mbF\ is \emph{normalized} if  $\{ 0,1\} \subseteq A$.

Let $A$ be a subset of $\mbF$ of cardinality $|A| \geq 2$.
Let $a_0, a_1 \in A$ with $a_0 \neq a_1$.  For all $a \in A$, let 
\[
a' = \frac{a-a_0}{a_1-a_0}
\]
and 
\[
A' = \left\{ \frac{a-a_0}{a_1-a_0}:a\in A \right\}.
\]
We have 
\[
a = (a_1 - a_0) a' + a_0 
\]
and 
\[
A = (a_1-a_0)\ast A' + a_0.
\] 
Note that $a'_0 = 0$ and $a'_1 = 1$, and so 
$\{ 0,1\} \subseteq A'$, that is, $A'$ is normalized. 
For all $(a_1,\ldots, a_h) \in A^h$ we have
\begin{align*}
\varphi(a_1,\ldots, a_h) 
& = \sum_{i=1}^h c_i  a_i  = \sum_{i=1}^h c_i \left( (a_i -a_0)a'_i  + a_0 \right) \\
& =   (a_1 - a_0) \sum_{i=1}^h c_i a'_i       + a_0 \sum_{i=1}^h c_i \\
& =    (a_1 - a_0) \varphi\left( a'_1 , \ldots,a'_h  \right) + a_0 \sum_{i=1}^h c_i .
\end{align*}
 Thus,
 \[
 \varphi(A) = (a_1-a_0)\ast \varphi(A') + a_0 \sum_{i=1}^h c_i.  
 \]
Because $a_1-a_0 \neq 0$, the set $A$ is a $\varphi$-Sidon set if and only if $A'$ is a $\varphi$-Sidon set. 
Thus, it suffices to consider only normalized subsets of \mbF.


\section{Binary linear forms}
Consider monic binary linear forms  $\varphi = x_1 + cx_2$ with $c \in \mbF^{\times}$.  
For every nonempty subset $A$ of $\mbF$, let 
\[                                         
\mce(A) = \left\{ c \in  \mbF^{\times} : \text{$A$ is not a $\varphi$-Sidon set for 
the  form  $\varphi = x_1 + cx_2$}  
\right\}.
\]
We shall prove that $\mce(A)$ is finite if $A$ is finite.  
Equivalently, every nonempty finite subset  $A$ of \mbF\ is a $\varphi$-Sidon set 
for all but finitely many monic binary linear forms $\varphi$.  
Moreover, we compute $\mce(A)$ explicitly for sets with at most 4 elements.  

Let $A$ be a nonempty subset of \mbF.   
Define the \emph{difference set}  
\[
D(A) = A-A = \{a' -a:a, a' \in A\}
\]
and the \emph{set of difference quotients}  
\[
D^{\ast}(A) = \left\{ \frac{d'}{d} : d,d' \in D(A) \text{ and } d\neq 0\right\}.
\]

\bt
Let $A$ be a nonempty subset of the field $\mbF$.  Then 
\[
\mce(A) \subseteq  D^{\ast}(A).
\]
\et

\begin{proof}
Let $c \in \mbF^{\times}$ and $\varphi = x_1+cx_2$.   If $c \in \mce(A)$,
then there exist $(a_1,a_2), (a'_1,a'_2)  \in A^2$
with $(a_1,a_2) \neq (a'_1,a'_2)$ such that 
\[
a_1 + ca_ 2 = \varphi (a_1,a_2) = \varphi (a'_1,a'_2) = a'_1 + ca'_ 2. 
\]
Equivalently, 
\[
c(a'_2 - a_2) = a_1-a'_1. 
\]
If $a_2 = a'_2$, then $a_1 = a'_1$ and  $(a_1,a_2) = (a'_1,a'_2)$, which is absurd. 
Therefore,  $a_2 \neq a'_2$ and 
\[
c = \frac{ a_1-a'_1}{a'_2 - a_2} \in D^{\ast}(A).
\]
This completes the proof.  
\end{proof}

\bc
Every nonempty finite subset of a field is a $\varphi$-Sidon set for all 
but finitely many monic binary forms.  
\ec

\begin{proof}
If $A$ is finite, then the sets $D(A)$ and $D^{\ast}(A)$ are also finite, 
and so $\mce(A)$ is finite.   
\end{proof}

We shall explicitly compute $\mce(A)$ for sets with 2, 3, or 4 elements.

\bl                                          \label{SidonInverse:theorem:AB}
If $A$ and $B$ are nonempty subsets of a field and $A\subseteq B$, then $\mce(A) \subseteq \mce(B)$. 
\el

\begin{proof}
Let $c \in \mce(A)$.  There exist distinct pairs $(a_1,a_2),(a'_1,a'_2) \in A^2$ 
such that $a_1+ca_2 = a'_1 + ca'_2$.  If $A \subseteq B$, then $(a_1,a_2),(a'_1,a'_2)$ 
are also distinct pairs in $B^2$ and so $c \in \mce(B)$.  This completes the proof.  
\end{proof}

For $c \in \mbF^{\times}$, let  
\[                      
E(c) = \left\{ \pm c,  \pm \frac{1}{c}  \right\}.  
\]

\bl                               \label{SidonInverse:lemma:Ec}
Let $A$ be a nonempty subset of the field $\mbF$. 
If $c \in \mce(A)$, then $E(c) \subseteq \mce(A)$.   
\el

\begin{proof} 
We associate to every $c \in \mbF^{\times}$ the monic binary linear forms 
\begin{align*} 
\varphi_1 & = x_1 + cx_2 \\
\varphi_2 & =  x_1 - cx_2 \\
\varphi_3 & =  x_1 + \frac{1}{c} x_2 \\
\varphi_4 & =  x_1 - \frac{1}{c} x_2.
\end{align*}

Let $(a_1,a_2), (a'_1,a'_2) \in A^2$.  We have 
\[
\varphi_1(a_1,a_2) = a_1+ca_2 = a'_1+ ca'_2  = \varphi_1(a'_1,a'_2) 
\]
if and only if 
\[
\varphi_2(a_1,a_2') = a_1 - ca'_2 = a'_1 - ca_2 = \varphi_2(a'_1,a_2) 
\]
if and only if 
\[
\varphi_3(a_2,a_1) = a_2 + \frac{1}{c} a_1  = a'_2 + \frac{1}{c} a'_1 = \varphi_3(a'_2,a'_1) 
\] 
if and only if 
\[
\varphi_4(a_2,a'_1) = a_2 - \frac{1}{c} a'_1  = a'_2 - \frac{1}{c} a_1 = \varphi_4(a'_2,a_1).
\]
Also, 
\[
(a_1,a_2) \neq (a'_1,a'_2)
\]
if and only if 
\[
(a_1,a'_2) \neq (a'_1,a_2)
\]
if and only if 
\[
(a_2,a_1) \neq (a'_2,a'_1)  
\]
if and only if 
\[
(a_2,a'_1) \neq (a'_2,a_1).   
\]
Thus, the set $A$ is a $\varphi$-Sidon set for some linear form 
$\varphi \in \{ \varphi_1, \varphi_2, \varphi_3, \varphi_4 \}$
 if and only if $A$ is a $\varphi$-Sidon set for all 
 $\varphi \in \{ \varphi_1, \varphi_2, \varphi_3, \varphi_4 \}$.  
Equivalently, $c \in \mce(A)$ implies $E(c) \subseteq \mce(A)$. 
This completes the proof.  
\end{proof}

\bt                                           \label{SidonInverse:theorem:A2}
Let $\mbF$ be a field and let $A = \{0,1\}$ be the normalized subset of $\mbF$ of cardinality 2.  
Then 
\[
\mce(0,1) = \{0, \pm 1\} = \{ 0 \} \cup E(1).
\]
\et

\begin{proof} 
For $c \in \mbF$ and $\varphi = x_1 + cx_2$, we have 
$\varphi(A) = \{0, 1 , c, 1+c\}$, and $|\varphi(A)| = 4$ if and only if $c \neq 0, \pm 1$.
Therefore, $\mce(A) = \{0, \pm 1\}$.  This completes the proof.  
\end{proof}

\bt                                          \label{SidonInverse:theorem:A3} 
Let $\mbF$ be a field and let $A = \{0,1,a\}$ be a normalized subset of $\mbF$ of cardinality 3.   
Then  
\begin{align*}
\mce(0,1,a) & =  \{  0 , \pm 1 \} \cup  
\left\{ \pm a,  \pm \frac{1}{a}  \right\}, \left\{ \pm (a-1),  \pm \frac{1}{a-1} \right\},
 \left\{  \pm \frac{a-1}{a}, \pm \frac{a}{a-1}  \right\} \\ 
& =  \mce(0,1)  \cup  E(a) \cup E(a-1) \cup E\left( \frac{a-1}{a} \right) 
\end{align*}
\et

\begin{proof} 
For the monic binary form $\varphi = x_1 + cx_2$,   
the set  $A = \{0,1,a\}$ is a $\varphi$-Sidon set if and only if the set 
\[
\varphi(A) = \{0,1,c,1+c,   a, ac,   a+c, 1+ac, a+ac  \}
\]
has cardinality 9.  There are $\binom{9}{2} = 36$ pairs of distinct elements  
of $\varphi(A)$.  
We obtain the set $\mce(A)$ by equating pairs of distinct elements 
of $\varphi(A)$ and solving for $c$.  
For example, for the elements $c \in \varphi(A)$ and $1+ac \in \varphi(A)$, we have 
$c=1+ac $ if and only if $c = -1/(a-1)$, and so  $c = -1/(a-1) \in \mce(A)$.  
Lemma~\ref{SidonInverse:lemma:Ec} implies that $E(a-1) \subseteq \mce(A)$.  
Computing $c$ for all 36 pairs determines $\mce(A)$.  
This completes the proof. 
\end{proof}

Note that the set 
\begin{align*} 
E(a) & \cup  E(a-1) \cup E\left( \frac{a-1}{a} \right) =  \\
&  \left\{ \pm a,  \pm \frac{1}{a}, \pm (a-1), \pm \frac{1}{a-1}, \pm \frac{a-1}{a}, \pm \frac{a}{a-1} \right\} 
\end{align*} 
is the set of linear fractional transformations associated 
with the following set of  matrices in $SL_2(\Z)$: 
\[
\pm \left\{  \bmat 1 & 0 \\ 0 & 1 \emat,  \bmat  0 & -1 \\ 1 & 0    \emat,  
  \bmat  1 & -1 \\ 0 & 1    \emat,   \bmat 0 & -1 \\1 & -1  \emat, 
     \bmat  1 & -1 \\ 1 & 0   \emat,   \bmat  -1 & 0 \\ 1 & -1   \emat 
\right\}.
\]

\bt                                        \label{SidonInverse:theorem:A4}
Let $\mbF$ be a field and let $A = \{0,1,a, b\}$ be a normalized subset of $\mbF$ of cardinality 4.   
Then 
\begin{align*}
\mce(0,1,a,b)  =  & \mce(0,1,a)  \cup \mce(0,1,b)   \\
&  \cup  E(a -b)  \cup E\left( \frac{a-b}{a} \right)   \cup E\left( \frac{a-b}{b} \right)  \\
&  \cup E\left( \frac{a-1}{b} \right) \cup E\left( \frac{b-1}{a} \right)  \cup E\left( \frac{a-1}{b-1} \right) \\
& \cup E\left(   \frac{a - 1}{a - b}  \right)   \cup E\left(   \frac{b - 1}{a - b}  \right)  \cup E\left( \frac{a}{b} \right). 
\end{align*}
\et

\begin{proof}
The set $A$ is a $\varphi$-Sidon set if and only if the set 
\begin{align*}
\varphi(A) = & \{0, 1, a, b, c, ca, cb, 1 + c, a + c, b + c,  \\ 
& \qquad   ac + 1,  ac + a, bc + 1, bc + b, bc + a, ac + b \}
\end{align*}
has cardinality 16.
As in the proof of Theorem~\ref{SidonInverse:theorem:A3}, we determine 
the elements $c \in \mce(A)$ by solving the $\binom{16}{2} = 120$ equations for $c$.  
Maple does the calculation.  Note that Lemma~\ref{SidonInverse:theorem:AB} 
implies that $\mce(0,1,a)  \cup \mce(0,1,b)  \subseteq \mce(0,1,a,b)$. 
\end{proof}

\section{Sidon sets for linear forms in $h \geq 2$ variables}

Here are two simple constructions of $h$-ary  linear forms $\varphi$ 
for which a finite set $A$ is a $\varphi$-Sidon set.

\bt
Let $A$ be a  finite set of complex numbers with $|A| \geq 2$.  Let 
\[
\delta = \min\{ |a - a'|: a, a' \in A \text{ and } a \neq a' \} 
\]
and
\[
\Delta = \max\{ |a - a'|: a, a' \in A \text{ and } a \neq a' \}. 
\]
Let $(c_i)_{i=1}^{\infty}$ be sequence of complex numbers 
such that 
\[
|c_j| > \frac{\Delta}{\delta}  \sum_{i=1}^{j-1} |c_i |
\]
for all $j \geq 2$.  The set $A$ is a $\varphi_h$-Sidon set for the 
$h$-ary linear form 
\[
\varphi_h = \sum_{i=1}^h c_i x_i 
\]  
for all $h \geq 2$. 
\et

\begin{proof}
If $A$ is not  a $\varphi_h$-Sidon set, then there exist distinct $h$-tuples 
 $(a_1,\ldots, a_h)$ and $(a'_1,\ldots, a'_h)$ in $A^h$ such that 
\[
\varphi_h(a_1,\ldots, a_h) = \varphi_h(a'_1,\ldots, a'_h).
\]
There is a largest integer $j \in \{1,\ldots, h\}$ 
such that $a_{j} \neq a'_{j}$.  It follows that 
\[
\sum_{i=1}^{j} c_i a_i = \sum_{i=1}^{j} c_i a'_i 
\]
and so 
\[
c_{j} \left( a'_{j} - a_{j} \right) = \sum_{i=1}^{j-1} c_i (a_i - a'_i ). 
\]
Therefore, 
\[
|c_j| \delta \leq \left| c_{j} \right| \left| a'_{j} - a_{j} \right| 
 \leq  \sum_{i=1}^{j-1} \left|  c_i \right|   \left| a_i - a'_i \right|  
 \leq \Delta   \sum_{i=1}^{j-1} \left|  c_i \right| < \delta  |c_j|
\]
which is absurd.  
This completes the proof.  
\end{proof}

The following construction appeared in~\cite{nath07m}.

\bt
Let $A$ be a nonempty finite set of positive integers and let 
$a^*= \max(A)$.    
For all integers $g > a^*$ and for all $h \geq 2$, 
the set $A$ is a $\varphi$-Sidon set for the $h$-ary linear form
\[
\varphi = x_1+gx_2+g^2x_3 + \cdots + g^{h-1} x_h.
\]
\et

\begin{proof}
This follows immediately from the uniqueness of the $g$-adic representation of a positive integer.  
\end{proof}

\section{Open problems} 
\benum

\item
Let $\Phi_h$ be the set of all monic $h$-ary linear forms 
$\varphi = x_1 + cx_2 + \cdots + c_h x_h$    with 
nonzero coefficients $c_i$ in the field \mbF. 
For every nonempty subset $A$ of $\mbF$, let 
\[                                         
\mce_h(A) = \left\{ \varphi \in \Phi_h : \text{$A$ is not a $\varphi$-Sidon set for 
some $h$-ary form  $\varphi$}  
\right\}.
\]
Compute the sets $\mce_h(A)$ for small subsets $A$ of \mbF.

\item
The \emph{representation function} of the set $A$ with respect to an $h$-ary 
linear form $\varphi$ with coefficients in the field \mbF\  is
\[
r_{A,\varphi}(b) = \card\left\{  (a_1,\ldots, a_h) \in A^h : \varphi(a_1,\ldots, a_h) = b \right\} 
\]
for $b \in \mbF$.  
Let $g \geq 1$.  The set $A$ is a \emph{$\varphi$-Sidon set  of order $g$}
 if $r_{A,\varphi}(b) \leq g$ for all $b \in \mbF$.  
The set $A$ is a $\varphi$-Sidon set if $A$ is a $\varphi$-Sidon set  of order $1$.  
It is of interest to investigate inverse problems for $\varphi$-Sidon sets of orders  $g \geq 2$.

\item
Consider the relationship between the exceptions sets $\mce_h(A)$ and subsets of the 
special linear group $SL_h(\mbF)$  both for $h =2$ and for $h \geq 3$.

\item
Consider the inverse Sidon problem for infinite sets.  Let $h \geq 2$.  
Does there exist an infinite set $A$ 
such that $A$ is not a $\varphi$-Sidon set for only finitely many $h$-ary linear forms $\varphi$?

\eenum

\def\cprime{$'$}
\providecommand{\bysame}{\leavevmode\hbox to3em{\hrulefill}\thinspace}
\providecommand{\MR}{\relax\ifhmode\unskip\space\fi MR }
\providecommand{\MRhref}[2]{%
  \href{http://www.ams.org/mathscinet-getitem?mr=#1}{#2}
}
\providecommand{\href}[2]{#2}


\begin{thebibliography}{1}

\bibitem{nath07m}
M.~B. Nathanson, \emph{Representation functions of bases for binary linear
  forms}, Funct. Approx. Comment. Math. \textbf{37} (2007), 341--350.

\bibitem{nath2021-200}
\bysame, \emph{Sidon sets for linear forms}, J. Number Theory (2022), 
{arXiv:2101:01034}, 2021.

\end{thebibliography}
\end{document}